\documentclass[12pt]{article}
\usepackage{amsthm,amsfonts,hyperref}
\usepackage{fullpage}

\newtheorem{thm}{Theorem}

\newtheorem{lem}{Lemma}

\newtheorem{rem}{Remark}
\newtheorem{dfn}{Definition}

\newcommand{\rank}{\mathrm{rank}}

\newcommand{\C}{\mathbb{C}}

\title{Spectral lower bounds for the quantum chromatic number of a graph}

\author{Pawel Wocjan\thanks{wocjan@cs.ucf.edu, Department of Computer Science, University of Central Florida, USA}\quad\quad Clive Elphick\thanks{clive.elphick@gmail.com, School of Mathematics, University of Birmingham, Birmingham, UK}}

\begin{document}

\maketitle

\begin{abstract}
The quantum chromatic number, $\chi_q(G)$,  of a graph $G$ was originally defined as the minimal number of colors necessary in a quantum protocol in which two provers that cannot communicate with each other but share an entangled state can convince an interrogator with certainty that they have a coloring of the graph. We use an equivalent purely combinatorial definition of $\chi_q(G)$ to prove that many spectral lower bounds for the chromatic number, $\chi(G)$, are also lower bounds for $\chi_q(G)$. This is achieved using techniques from linear algebra called pinching and twirling. We illustrate our results with some examples.

\end{abstract}

\section{Introduction}

For any graph $G$ let $V$ denote the set of vertices where $|V| = n$, $E$ denote the set of edges where $|E| = m$, $A$ denote the adjacency matrix, $\chi(G)$ denote the chromatic number and $\omega(G)$ the clique number. Let $\mu_1 \ge \mu_2 \ge ... \ge \mu_n$ denote the eigenvalues of $A$, and then the inertia of $G$ is the ordered triple $(n^+, n^0, n^-)$, where $n^+, n^0$ and $n^-$ are the numbers of positive, zero and negative eigenvalues of $A$, including multiplicities. Note that $\rank(A) = n^+ + n^-$ and $\mathrm{nullity}(A) = n^0$. Let $s^+$ and $s^-$ denote the sum of the squares of the positive and negative eigenvalues of $A$, respectively.

Let $D$ be the diagonal matrix of vertex degrees, and let $L = D - A$ denote the Laplacian of $G$ and $Q = D  + A$ denote the signless Laplacian of $G$. The eigenvalues of $L$ are $\theta_1 \ge ... \ge \theta_n = 0$ and the eigenvalues of $Q$ are $\delta_1 \ge ... \ge \delta_n$.

Let $\chi_q(G)$ and $\chi_q^{(r)}(G)$ denote the quantum and rank-$r$ quantum chromatic numbers, as defined by Cameron \emph{et al} \cite{cameron07}, where $\chi_q(G) = \min_r{(\chi_q^{(r)}(G))}$. It is evident that $\chi_q(G) \le \chi(G)$, and Cameron \emph{et al} \cite{cameron07} exhibit a graph on 18 vertices and 44 edges with chromatic number 5 and quantum chromatic number 4. Mancinska and Roberson \cite{mancinska162} have subsequently found a graph on 14 vertices with $\chi(G) > \chi_q(G)$, and they suspect this is the smallest possible example.
 
It is helpful to have a purely combinatorial definition of  the quantum chromatic number, and the following definition is due to \cite[Definition 1]{mancinska162}.

For a positive integer $c$, let $[c]$ denote the set $\{0,1,\ldots,c-1\}$. For $d>0$, let $I_d$ and $0_d$ denote the identity and zero matrices in $\C^{d\times d}$.
 
\begin{dfn}
A quantum $c$-coloring of the graph $G=(V,E)$  is a collection of orthogonal projectors $\{ P_{v,k} : v\in V, k\in [c]\}$ in 
$\C^{d\times d}$ such that
\begin{itemize}
\item for all vertices $v\in V$
\begin{eqnarray}
\sum_{k\in[c]} P_{v,k} & = & I_d \quad\quad \mathrm{(completeness)} \label{eq:complete}
\end{eqnarray}
\item for all edges $vw\in E$ and for all $k\in[c]$
\begin{eqnarray}
P_{v,k} P_{w,k}           & = & 0_d \quad\quad \mathrm{(orthogonality)} \label{eq:orthogonal}
\end{eqnarray}
\end{itemize}
The quantum chromatic number $\chi_q(G)$ is the smallest $c$ for which the graph $G$ admits a quantum $c$-coloring for some dimension $d>0$.
\end{dfn}

According to the above definition, any classical $c$-coloring can be viewed as a $1$-dimensional quantum coloring, where we set $P_{v,k}=1$ if vertex $v$ has color $k$ and
we set $P_{v,k}=0$, otherwise. Therefore, quantum coloring is a relaxation of  classical coloring. As noted in \cite{mancinska162}, it is surprising that the quantum chromatic number can be strictly and even exponentially smaller than the chromatic number for certain families of graphs.

\section{Spectral lower bounds for the chromatic number}
 
Most of the known spectral lower bounds for the chromatic number can be summarized as follows:
 
\begin{equation}\label{bounds}
1 + \max\left(\frac{\mu_1}{|\mu_n|} , \frac{2m}{2m- n\delta_n} , \frac{\mu_1}{\mu_1 - \delta_1 + \theta_1} , \frac{n^\pm}{n^\mp} ,  \frac{s^\pm}{s^\mp}\right) \le \chi(G) , 
\end{equation}
where, reading from left to right, these bounds are due to Hoffman \cite{hoffman70}, Lima \emph{et al} \cite{lima11}, Kolotilina \cite{kolotilina11}, Elphick and Wocjan \cite{elphick17}, and Ando and Lin \cite{ando15}.  It should be noted that Nikiforov \cite{nikiforov07} pioneered the use of non-adjacency matrix eigenvalues to bound $\chi(G)$.

Let $c$ denote the number of colors used in a (classical) coloring. The authors (\cite{wocjan13} \cite{elphick15}, \cite{elphick17}), provided proofs of all of the bounds in (\ref{bounds}) except the last one using the following equality:

\begin{equation}\label{conversion}
\sum_{\ell\in[c]} U^\ell A (U^{\dagger})^\ell = 0_n,
\end{equation}
where $U$ is a diagonal unitary matrix in $\C^{n\times n}$ whose entries are $\chi$th roots of unity and $\dagger$ denotes the conjugate transpose.  The last bound in (\ref{bounds}) was proved in \cite{ando15} essentially using the following equality:
\begin{equation}\label{partition}
\sum_{k\in[c]} P_k A P_k = 0_n\,, 
\end{equation}
where $P_k$ are orthogonal projectors in $\C^{n\times n}$ that are diagonal in the \emph{standard basis} and their sum  $\sum_{k\in[c]}P_k$ is equal to the identity matrix $I_n$. 

A quantum $c$-coloring is an operator relaxation of a classical $c$-coloring. The latter corresponds to the special case when the dimension $d$ of the relevant Hilbert space is $1$.
We will show that the existence of a quantum $c$-coloring in dimension $d$ implies the existence of suitable orthogonal projectors $P_k$ and a suitable unitary matrix $U$ in $\C^{n\times n}\otimes\C^{d\times d}$ such that the above equalities hold for $A\otimes I_d$. Once these equalities are established, we can use the same approaches as in the above papers to prove that all bounds in (\ref{bounds}) are also lower bounds for the quantum chromatic number.

We note that all bounds are also valid for weighted adjacency matrices of the form $W \circ A$, where $W$ is an arbitrary Hermitian matrix and $\circ$ denotes the Hadamard product (also called the Schur product).

\section{Pinching and twirling}

We start by defining two operations from linear algebra: pinching and twirling.

\begin{rem}
The following observation is fairly obvious but important. Let $\{Q_k : k \in[c]\}$ be any collection of orthogonal projectors in $\mathbb{C}^{m\times m}$ that form a resolution of the identity matrix, that is, 
\[
\sum_{k\in[c]} Q_k = I_m\,.
\] 
Then, the orthogonal projectors are necessarily mutually orthogonal, that is, $Q_k Q_\ell = 0$ for $k,\ell\in[c]$ with $k\neq \ell$.
\end{rem}

The following definition of pinching can be found in \cite[Problem II.5.5.]{bhatiaBook}.

\begin{dfn}[Pinching]\label{dfn:pinching}
Let $\{Q_k : k \in [c]\}$ be any collection of orthogonal projectors in $\C^{m\times m}$ that form a resolution of the identity matrix.
Then, the operation $\mathcal{C}$ that maps an arbitrary matrix $X\in\C^{m \times m}$ to
\[
\mathcal{C}(X) = \sum_{k\in[c]} Q_k X Q_k
\]
is called pinching. We say that the pinching $\mathcal{C}$ annihilates $X$ if $\mathcal{C}(X)=0_m$. 
\end{dfn}

Let $\{e_i : i\in [m]\}$ denote the standard basis of $\C^m$.  The basis vector $e_i$ has $1$ in the $i$th position and $0$ in all other positions.

\begin{rem}[Partitioning and pinching]\label{partition}
Assume that the row and column indices of matrices $X\in\C^{m\times m}$ are partitioned into the following $c$ nonempty sets $S_k=\{s_k,\ldots,s_{k+1}-1\}$ for $k\in[c]$ for
given $0=s_0<s_1<\ldots<s_{c-1}<s_{c}=m$. Let
\[
X=\left(
\begin{array}{cccc}
X_{0,0} & X_{0,1} & \cdots & X_{0,c-1} \\
X_{1,0} & X_{1,1} & \cdots & X_{1,c-1} \\
\vdots    & \vdots    & \ddots & \vdots \\
X_{c-1,0} & X_{c-1,1} & \cdots & X_{c-1,c-1}
\end{array}
\right)
\]
be the corresponding partitioned matrix.

Let $\{ P_k : k \in [c]\}$ be a collection of projectors in $\C^{m\times m}$, where $P_k$ denote the projectors onto the subspaces
\[
\mathrm{span}\{ e_i : i \in S_k \}.
\] 
Then, $P_k X P_\ell$ correspond to the submatrices $X_{k,\ell}$ of  the above partition of $X$. Let $\mathcal{C}$ be the pinching
defined by the collection $\{ P_k : k \in [c] \}$ of the above orthogonal projectors. Then, 
\[
\mathcal{C}(X) = 
\left(
\begin{array}{cccc}
X_{0,0} & 0           & \cdots & 0 \\
0           & X_{1,1} & \cdots & 0 \\
\vdots    & \vdots    & \ddots & \vdots \\
0 & 0     & \cdots    & X_{c-1,c-1}
\end{array}
\right)\,.
\]
\end{rem}

\begin{dfn}[Twirling]
Let $\{U_\ell : \ell \in [c]\}$ be a collection of arbitrary unitary matrices in $\C^{m\times m}$. Borrowing terminology from quantum information theory,
we call the operation $\mathcal{D}$ that maps an arbitrary matrix $X\in\C^{m\times m}$ to
\[
\mathcal{D}(X) = \frac{1}{c} \sum_{\ell\in [c]} U_\ell X U_\ell^\dagger
\]  
twirling. We say that the twirling $\mathcal{D}$ annihilates $X$ if $\mathcal{D}(X)=0_m$.
\end{dfn}

It was shown in \cite{bhatia} that twirling can be constructed from pinching in a straightforward way so that both operations have the same effect. Observe that in this construction the unitary matrices $U_\ell$ above can be chosen to be powers of one unitary matrix $U$, that is, we have $U_\ell=U^\ell$.  The special case when there are only two projectors is mentioned in \cite[Problem II.5.4]{bhatiaBook}.

\begin{lem}\label{lem:twirling}
It is known that pinching $\mathcal{C}$ defined in Definition~\ref{dfn:pinching} can also be realized as twirling $\mathcal{D}$ as follows. 
Let $\omega=e^{2\pi i/c}$ be a $c$th root of unity and 
\[
U = \sum_{k\in [c]} \omega^{k} Q_k\,.
\]
Then, twirling defined by
\begin{equation}\label{eq:averaging}
{\mathcal D}(X) = \frac{1}{c} \sum_{\ell\in[c]} U^\ell X (U^{\dagger})^\ell\,. 
\end{equation}
satisfies
\[
{\mathcal C}(X) = {\mathcal D}(X)
\]
for all matrices $X\in\C^{m\times m}$. 
\end{lem}
\proof{The $\ell$th power of $U$ is equal to 
\[
U^\ell = \sum_{k\in [c]} \omega^{k\cdot \ell} Q_k
\]
since the projectors $Q_k$ are mutually orthogonal to each other. We obtain
\begin{eqnarray*}
{\mathcal D}(X) 
& = &
\frac{1}{c} \sum_{\ell\in[c]} U^\ell X (U^{\dagger})^\ell \\
& = & 
\frac{1}{c} \sum_{k,k'\in[c]} \sum_{\ell\in[c]} \omega^{(k-k')\cdot \ell} P_k X P_{k'} \\
& = & 
\sum_{k\in[c]} P_k X P_k = {\mathcal C}(X).
\end{eqnarray*}
In the last step, we use that $\sum_{\ell\in[c]} \omega^{(k-k')\cdot \ell}=c\cdot\delta_{k,k'}$, where $\delta_{k,k'}$ denotes the Kronecker delta. 
}

\section{Pinching from quantum coloring}

We will show how to construct pinchings from quantum colorings. In particular, we will show that if there exists a quantum $c$-coloring 
in dimension $d$, then there exists a pinching with $c$ orthogonal projectors that annihilates $A\otimes I_d$.

Let $\{e_v : v \in V\}$ denote the standard basis in $\C^n$, where $n=|V|$.
Denote the entries of $A$ by $a_{uv}$, where $u,v\in V$ enumerate the rows and columns, respectively. We have 
\[
A = \sum_{v,w\in V} a_{vw} e_v e_w^\dagger\,,
\]
where $a_{vw}=e_v^\dagger A e_w$.

\begin{thm}\label{thm:main}
Let $\{ P_{v,k} : v\in V, k\in [c]\}$ be an arbitrary quantum $c$-coloring of $G$ in $\C^d$. Then, the 
following block-diagonal orthogonal projectors 
\[
P_k = \sum_{v\in V} e_v e_v^\dagger \otimes P_{v,k} \in \C^{n\times n}\otimes\C^{d\times d} 
\]
form a resolution of the identity matrix. Moreover, the corresponding pinching operation $\mathcal{C}$ 
\begin{itemize}
\item annihilates $A\otimes I_d$, and 
\item leaves $E\otimes I_d$ invariant for all diagonal matrices $E\in\C^{n\times n}$.
\end{itemize}
\end{thm}

\proof{They form a resolution of the identity matrix because
\[
\sum_{k\in[c]} P_k  = \sum_{v\in V} e_v e_v^\dagger \otimes \sum_{k\in [c]} P_{v,k} = \sum_{v\in V} e_v e_v^\dagger \otimes I_d = I_n \otimes I_d\,,
\]
where we make use of the completeness condition (1) that $\sum_{k\in [c]} P_{v,k} = I_d$ for all vertices $v\in V$.

The corresponding pinching operation $\mathcal{C}$ annihilates $A \otimes I_d$ because
\begin{eqnarray*}
& &
\mathcal{C}(A\otimes I_d) \\ 
& = &
\sum_{k\in[c]} P_k (A \otimes I_d) P_k \\
& = & 
\sum_{k\in[c]} \left(\sum_{v\in V} e_v e_v^\dagger \otimes P_{v,k}\right) (A \otimes I_d) \left( \sum_{w\in V} e_w e_w^\dagger \otimes P_{w,k} \right) \\
& = &
\sum_{k\in[c]} \sum_{v,w\in V} a_{vw} \cdot e_v e_w^\dagger \otimes P_{v,k} P_{w,k} \\
& = &
\sum_{k\in[c]} \left( \sum_{vw\in E} 1 \cdot e_v e_w^\dagger \otimes 0_d + \sum_{vw\not\in E} 0 \cdot e_v e_w^\dagger \otimes P_{v,k} P_{w,k} \right) \\
& = & 0\,,
\end{eqnarray*}
where we made use of the orthogonality condition (2) $P_{v,k} P_{w,k} = 0$ for all $vw\in E$ (or equivalently, for all pairs $v,w\in V$ with $a_{vw}=1$).

The property that $\mathcal{C}$ leaves $E \otimes I_d$ invariant for all diagonal matrices $E\in\C^{d\times d}$ is verified similarily.
}

%\begin{lem}\label{lem:invariance}
%Let $G=(V,E)$ be an arbitrary graph that has a quantum $c$-coloring in dimension $d$.  Let $U$ be the unitary matrix obtained by (i) using Theorem~\ref{thm:main} to construct a pinching from the quantum coloring and (ii) Lemma~\ref{lem:twirling} to translate this pinching into a twirling. Then, the unitary matrix $U$ satisfies the following equation: 
%\[
%U^\ell (D \otimes I_d) (U^\dagger)^\ell = D\otimes I_d\,,
%\]
%for arbitrary diagonal matrices 
%\[
%D=\sum_{v\in V} d_v \cdot e_v^\dagger e_v 
%\]
%in $\C^{n\times n}$ and all $\ell\in[c]$. 
%\end{lem}

%\begin{proof}
%The proof of this Lemma follows very closely that of Theorem~\ref{thm:main}.
%\end{proof}

Theorem~\ref{thm:main} shows that the existence of a quantum $c$-coloring of a graph $G=(V,E)$ with adjacency matrix $A$ in dimension $d$ implies the existence of a pinching operation $\mathcal{C}$ that annihilates $A\otimes I_d$ and leaves $E\otimes I_d$ invariant for all diagonal matrices $E$.  The converse direction is shown in the remark below \cite{robertson}. 

\begin{rem}\label{rem:converse}
Theorem 4.25 in  \cite{watrous} shows that the fixed points of any completely positive trace-preserving unital map commute with the Kraus operators of the map. In the present case, the completely positive trace-preserving map is the pinching $\mathcal{C}$, the Kraus operators are the orthogonal projectors $P_k$, and the fixed points of $\mathcal{C}$ are $E\otimes I_d$, where $E\in\C^{d\times d}$ is an arbitrary diagonal matrix. This result implies that the $P_k$ commute with $E \otimes I_d$, which in turn implies that the $P_k$ are block diagonal. These blocks are indexed by the vertices $v$ of $G$, so we can refer to them as $P_{v,k}$ for $v \in V$.

Using the block-diagonal nature of the $P_k$, it is now easy to show that these $P_{v,k}$ yield a quantum $c$-coloring of $G$: 
\begin{itemize} 
\item The property $\sum_{k\in[c]} P_k = I_{nd}$ implies $\sum_{k\in[c]} P_{v,k} = I_{d}$ for $v\in V$, thus giving the completeness condition in (\ref{eq:complete}).
\item The property $\mathcal{C}(A \otimes I) = 0_{nd}$ implies $\sum_{k\in[c]} P_{v,k} P_{w,k} = 0_d$ for all $vw\in E$.  Let $\ell\in[c]$ be fixed by arbitrary.  Multiplying the latter equation by $P_{v,\ell}$ from the left and by $P_{w,\ell}$ from the right shows that the summand $P_{v,\ell} P_{w,\ell}$ must be zero, thus giving the orthogonality condition in (\ref{eq:orthogonal}).
\end{itemize}

The pinching operation described in this paper can therefore be regarded as an algebraic  reformulation of quantum coloring.
\end{rem}

\section{Lower bounds on quantum chromatic number}

Using Theorem~\ref{thm:main} and Lemma~\ref{lem:twirling}, it is possible to show that all the bounds in (\ref{bounds}) are also lower bounds on the quantum chromatic number.  

We demonstrate this explicitly for the bound
\[
1 +  \frac{2m}{2m- n\delta_n}\,,
\]
where $\delta_n$ is the minimum eigenvalue of the signless Laplacian $Q=D+A$.

Assume that there exists a quantum $c$-coloring in dimension $d$. Let $\{Q_k : k \in [c] \}$ denote projectors defining a pinching as in Theorem~\ref{thm:main} and $U^\ell=\sum_{k\in [c]} \omega^{k\cdot \ell} Q_k$ denote the corresponding twirling unitaries as defined in Lemma~\ref{lem:twirling}.   

The proof is almost identical to the proof for the chromatic number \cite{elphick15}.  
We use the identity $D-Q=-A$. We have:
\begin{eqnarray*}
A\otimes I_d 
& = &
\sum_{\ell=1}^{c-1} U^\ell (-A\otimes I_d) (U^{\dagger})^\ell \\
& = & 
\sum_{\ell=1}^{c-1} U^\ell \left( (D-Q)\otimes I_d \right) (U^{\dagger})^\ell \\
& = &
(c-1)(D\otimes I_d) - \sum_{\ell=1}^{c-1} U^\ell (Q\otimes I_d) (U^{\dagger})^\ell.
\end{eqnarray*}
Define the column vector $v=\frac{1}{\sqrt{nd}}(1,1,\ldots,1)^T$. Multiply the left and right most sides of the above matrix equation by $v^\dagger$ from the left and by $v$ from the right to
obtain
\[
\frac{2m}{n} = v^\dagger (A\otimes I_d)v = (c-1)\frac{2m}{n} - \sum_{\ell=1}^{c-1} v^\dagger U_\ell (Q \otimes I_d) U^\dagger_\ell v \le (c-1)\frac{2m}{n} - (c-1)\delta_n\,.
\]
This uses that $v^\dagger (A\otimes I_d) v = v^\dagger(D\otimes I_d) v = 2m/n$, which is equal to the sum of all entries of respectively $A$ and $D$ divided by $n$ due to the special form of $v$, and that $v^\dagger U_\ell (Q\otimes I_d) U^\dagger_\ell v \ge \lambda_{\min}(Q) = \delta_n$. 

The other bounds in (\ref{bounds}) can be  shown to be lower bounds for $\chi_q(G)$, by similarly modifying  the proofs of these bounds, so that they can be applied to $A\otimes I_d$ instead of $A$. 

\section{Implications for the quantum chromatic number}

We now discuss some implications of the bounds in (\ref{bounds}) being lower bounds on the quantum chromatic number. In particular, we discuss applications of the inertia bound: $\chi_q(G)\ge 1 + n^\pm/n^\mp$.

Let $\chi_v(G)$ denote the vector chromatic number of $G$ and $\theta(\overline{G})$ denote the Lov\'asz theta function of the complement of $G$. It is known that:
\[
1 + \frac{\mu_1}{|\mu_n|} \le \chi_v(G) \le \theta(\overline{G}) \le \chi_q(G)\,,
\]
where these inequalities (from left to right) are due to Bilu \cite{bilu06} , Karger \emph{et al} \cite{karger98} and Mancinska and Roberson \cite{mancinska16}. So it is already known that the Hoffman lower bound for $\chi(G)$ is also a lower bound for $\chi_q(G)$. Consequently it is the Lima \emph{et al}, Kolotilina, Ando and Lin and inertial bounds which are new. Experimentally the inertial bounds usually perform best in this context, so it is these bounds we focus on below. (The Lov\'asz theta bound is in general more difficult to compute than spectral bounds.)

In order for the inertial bounds  to reveal potentially new information about $\chi_q(G)$ it is necessary for:

\[
1 + \max\left(\left\lceil\frac{n^+}{n^-}\right\rceil , \left\lceil\frac{n^-}{n^+}\right\rceil\right) > \max{\left(\omega(G) , 1 + \left\lceil\frac{\mu_1}{|\mu_n|}\right\rceil\right)}.
\]

This is the case for many graphs, and we tabulate in Table 1 a few examples.

\begin{table}[ht]
\caption{Inertia vs Hoffman bounds}
\centering
\begin{tabular}{c c c c c c c}
\hline \hline
Graph & $n$ & $\chi$ & $\chi_q$ & Inertia & Hoffman & $\omega$\\[0.5ex]
\hline
Cyclotomic & $13$ & $4$ & $4$ & $3.25$ & $2.51$ & $2$\\
Clebsch & $16$ & $4$ & $4$ & $3.2$ & $2.67$ & $2$\\
Generalised Quadrangle(2,4) & $ 27 $ & $ 6 $ & $ \ge5 $ & $ 4.5 $ & $ 3 $ & $ 3 $\\
Non-Cayley Transitive(28,3) & $ 28 $ & $ 4 $ & $ 4 $ & $ 3.1 $ & $ 2.67 $ & $ 2 $\\
\hline
\end{tabular}
\end{table}

So, for example, the Hoffman/Bilu bound implies that the Clebsch graph has $\chi_q(G) \ge 3$ but the inertial bound implies $\chi_q(G) = 4$. More generally, $\chi_q(G) = \chi(G)$ if the ceiling of the inertia bound equals $\chi(G)$.

\section{Conclusion}
Our results can be interpreted as follows. They demonstrate that any existing or future general lower bound on the minimum number of operations required for pinching or twirling to annihilate a given matrix representation of a graph, becomes automatically a lower bound on the quantum chromatic number of that graph.

\section{Acknowledgement}

We would like to thank David Roberson for insightful comments on an earlier version of this paper, and in particular for Remark~\ref{rem:converse}.

This research has been partially supported by the National Science Foundation Award \#1525943 ``Is the Simulation of Quantum Many-Body Systems Feasible on the Cloud?''

\end{document}